\documentclass[12pt]{article}
\input isolatin1.sty
 \usepackage{a4}
 \usepackage{amsmath}
 \usepackage{amssymb}
\usepackage{amsfonts}
\usepackage{amsthm}
\usepackage{graphicx}
\newtheorem{thm}{Theorem}[section]
\newtheorem{cor}[thm]{Corollary}
\newtheorem{lem}[thm]{Lemma}

\newtheorem{rem}[thm]{Remark}

\renewcommand{\a}{\alpha}

\renewcommand{\t}{\theta}

\numberwithin{equation}{section}  

\title{Fibonacci numbers, Euler's 2-periodic continued fractions and moment
sequences
\footnote{This
  work was done during the visit of the second author to the University
of Copenhagen  partially supported by  D.G.E.S,
  ref. BFM2006-13000-C03-01, FQM-262, P06-FQM-01735, FQM-481 ({\em Junta de Andaluc\'{\i}a}), and by
  grant
  272-07-0321 from the Research Council for Nature and Universe.}}
\author{Christian Berg $^{\dagger}$ and Antonio J. Dur\'{a}n  $^{\ddagger}$\\
   \footnotesize $\dagger$ \ Institut for Matematiske Fag. K\o benhavns
   Universitet \\
    \footnotesize Universitetsparken 5; DK-2100 K\o benhavn \o,
   Denmark. berg@math.ku.dk \\  \footnotesize
   $\ddagger$ \  Departamento de An\'{a}lisis Matem\'{a}tico.
   Universidad de Sevilla \\
   \footnotesize Apdo (P. O. BOX) 1160; E-41080 Sevilla, Spain. duran@us.es \\
\ \ }

\begin{document}
\maketitle

\begin{abstract}
We prove that certain  sequences of finite continued fractions associated with a 2-periodic continued fraction with period $a,b>0$
are moment sequences of  discrete signed measures supported in the
interval $[-1,1]$, and we give necessary and sufficient conditions
in order that these measures are positive.
For $a=b=1$ this proves that the sequence of ratios
$F_{n+1}/F_{n+2},\,n\ge 0$ of  consecutive Fibonacci numbers is a
moment sequence.

\end{abstract}



\section{Introduction}

The last chapter of Euler's masterpiece {\sl Introductio in Analysin
Infinitorum} (vol. I), in English version \cite{E}, is devoted to continued fractions.
Euler considered there 1- and 2-periodic continued fractions
to get rational approximations to square roots of natural numbers. We
say that a continued fraction of the form  
\begin{equation}\label{pcf}
\frac{1}{a_1+\displaystyle \frac{1}{a_2+\displaystyle \frac{1}{a_3+\displaystyle \frac{1}{a_4+_{\ddots}}}}} \quad
\end{equation}
is $k$-periodic if $a_{j+k}=a_j$ for $j=1,2,\ldots$.

If a natural number $n$ is sum of two squares $n=m^2+l^2$ of natural
numbers, then for $a=2m/l$ the convergents of the 1-periodic
 continued fraction
\begin{equation}\label{1pcf}
\frac{1}{a+\displaystyle \frac{1}{a+\displaystyle \frac{1}{a+\displaystyle \frac{1}{a+_{\ddots}}}}} \quad
\end{equation}
give rational approximations for $\sqrt n$. Indeed, this continued fraction converges to the positive root of
$x^2+ax-1$:
$$
x=\frac{-a+\sqrt{a^2+4}}{2}=-\frac{m}{l}+\frac{\sqrt{n}}{l}.
$$
For example, for $n=5$, Euler took $m=1, l=2$, and so $a=1$; the continued fraction
(\ref{1pcf}) converges then to $(\sqrt{5}-1)/2$. In this example, the
rational approximations are
\begin{equation}\label{rafn}
\frac{0}{1},\quad \frac{1}{1}, \quad \frac{1}{2}, \quad \frac{2}{3},\quad \frac{3}{5}, \quad \frac{5}{8}, \quad \cdots
\end{equation}
which are ratios of consecutive Fibonacci numbers $F_n,n\ge 0$:
$F_0=0,F_1=1,F_{n+1}=F_n+F_{n-1},n\ge 1$. For a relation between  quotients
of consecutive Fibonacci numbers and electrical networks see
\cite[p.43]{Ko}.

In \cite{Is1} Ismail introduced a generalization $F_n(\t)$ of the
Fibonacci numbers as the solution to the difference equation
\begin{equation}\label{eq:genfib}
x_{n+1}=2\sinh\t\, x_n+x_{n-1},\quad n\ge 1,
\end{equation}
with the initial conditions $x_0=0,x_1=1$. Here $\t$ can be any
positive real number. The generalized Fibonacci
numbers are natural numbers, when $2\sinh\t$ is a natural number, and
the Fibonacci numbers correspond to $2\sinh\t=1$. 

When $n$ is not a sum of two squares, Euler considered a 2-periodic
continued fraction of the form
\begin{equation}\label{2pcf}
\frac{1}{a+\displaystyle \frac{1}{b+\displaystyle \frac{1}{a+\displaystyle \frac{1}{b+_{\ddots}}}}}
\end{equation}
to get rational approximations for $\sqrt n$. To do that, Euler  was implicitly using the Pell equation $m^2n=d^2-1$.
That equation was a proposal of Fermat to the British mathematicians. It was studied by Wallis and Brouncker. Euler thought
that Brouncker's  results, published by Wallis in his \textit{Algebra},
were due to Pell, and baptized it 
\textit{Pell's equation}.
Euler also studied Pell's equation, but it was Lagrange who solved it completely.
Assuming that $n$ is not a square of a natural number,
 the Pell equation $m^2n=d^2-1$
has always infinitely many solutions $m,d\in \mathbb N$.
 By taking positive integers $a,b$ for which
$ab=2d-2$, the 2-periodic continued fraction above gives rational approximations for $\sqrt n$. Indeed,
that continued fraction converges to the positive root of
$ax^2+abx-b$:
$$
x=\frac{-ab+\sqrt{a^2b^2+4ab}}{2a}=-\frac{d-1}{a}+\frac{m\sqrt{n}}{a}.
$$
For instance, for $n=7$, Euler took $m=3, d=8, a=2, b=7$, and the continued fraction (\ref{2pcf}) converges
then to $(-7+3\sqrt{7})/2$.

For a 2-periodic continued fraction  (\ref{2pcf}), where $a,b$ are positive real numbers, we associate the sequence
$(s_n(a,b,w))_n$ of finite continued fractions (to simplify the notation we remove the dependence on $a,b$ and $w$):
\begin{equation}\label{partfrac}
s_0=w, \quad s_1=\frac{1}{a+w},\quad s_2=\frac{1}{a+\displaystyle \frac{1}{b+w}},\quad 
s_3=\frac{1}{a+\displaystyle \frac{1}{b+\displaystyle \frac{1}{a+w}}}, \quad \cdots
\end{equation}
where $w\ge 0$. For $w=0$ we obtain the sequence of convergents to (\ref{2pcf}).
Since
$$
s_{n+2}=\frac{1}{a+\displaystyle \frac{1}{b+s_n}},
$$
we find
$$
s_n=\frac{N_n}{D_n},\quad n\ge 0,
$$
where the sequences $(N_n)_n$ and $(D_n)_n$ can be defined recursively:
\begin{eqnarray}\label{2suc}
N_{n+2}&=&bD_n+N_n,\quad n\ge 0,\\ \label{2suc1}
D_{n+2}&=&abD_n+aN_n+D_n, \quad n\ge 0,
\end{eqnarray}
with initial conditions $N_0=w$, $N_1=1$, $D_0=1$, $D_1=a+w$.

When $a=b$, by setting $D_{-1}=w$, we actually have $N_{n}=D_{n-1}$, $n\ge 0$, and then the sequence $(D_n)_n$ can be defined recursively
in the form
\begin{equation}\label{eq:KM}
D_{n+1}=aD_n+D_{n-1},\quad n\ge 0,
\end{equation}
with initial conditions $D_{-1}=w$, $D_0=1$. This is the same
difference equation as \eqref{eq:genfib}. The general difference
equation of second order with constant coefficients $x_{n+1}=ax_n+bx_{n-1}$ has been studied  by
Kalman and Mena in \cite{K:M}.

The moments $\mu_n$, $n\ge 0$, of a positive Borel measure $\mu$ on
the real line are defined by $\mu _n=\int _{\mathbb R} t^nd\mu(t) $, $n\ge 0$,
assuming that these integrals are finite. Sequences that are moments of positive measures on $[0,\infty)$
were characterized by Stieltjes in his fundamental memoir \cite{St}
by certain quadratic forms being non-negative. In \cite{Ham} Hamburger extended the results of
Stieltjes to moment sequences of positive measures concentrated on the whole real line.
Later Hausdorff, cf. \cite{Ha}, characterized moment sequences
of measures concentrated on the unit interval $[0,1]$ by
complete monotonicity, see also \cite{Ak}, \cite{W}, \cite{BCR}. Moment
sequences of positive measures on the interval $[-1,1]$ can be
characterized as bounded Hamburger moment sequences, i.e. bounded
sequences $(s_n)_n$ such that all the Hankel matrices
\begin{equation}\label{eq:Hankel}
\mathcal H_n=\left(s_{i+j}\right)_{i,j=0}^n,\quad n=0,1,\ldots
\end{equation} 
are positive semidefinite.

The purpose of this note is to characterize the
sequences (\ref{partfrac}) of finite continued fractions $(s_n)_n$,
 which are bounded moment sequences.

In the following $\delta_c$ denotes the Dirac measure having the mass 1 at the point $c\in\mathbb R$.

\begin{thm}\label{thm:1} For $a,b,w\in \mathbb R$, $a,b>0, w\ge 0$, define
\begin{equation}\label{exprqab}
q=\frac{2+ab-\sqrt{a^2b^2+4ab}}{2},\quad \alpha =\frac{q(aw-(1-q))}{qaw+1-q},\quad \beta =\frac{qa-(1-q)w}{a+(1-q)w}.
\end{equation}
Then $0<q<1,|\alpha|,|\beta|<1$ and $(s_n(a,b,w))_n$ is the sequence
of moments of the discrete signed measure $\rho$ supported in $[-1,1]$ and defined by
\begin{equation}\label{defro}
\rho=\frac{1}{a}(1-q)\delta_1+\frac{1}{2a}\left( \frac{1}{q}-q\right)\sum_{k=0}^\infty
\left( \alpha ^{k+1}+\beta^{k+1}\right)\delta_{q^{k+1}}+\left( \alpha ^{k+1}-\beta^{k+1}\right)\delta_{-q^{k+1}}.
\end{equation}
The measure $\rho $ is positive if and only if $a\ge b$,
$$
w\ge -b/2+\sqrt{(b/2)^2+b/a},\quad \mbox{and} \quad w\ge -(a+b)/4+\sqrt{((a+b)/4)^2+1}.
$$
\end{thm}

In particular for $a\ge b$ and $w=1$, the measure $\rho$ is always positive.
Taking $a=b=2\sinh\t>0$, and $w=1/a$ we get the following result:

\begin{cor}\label{thm:genfib} The sequences
\begin{equation}\label{eq:quofib}
\frac{F_{n+1}(\t)}{F_{n+2}(\t)},\quad
\frac{F_{n+3}(\t)}{F_{n+2}(\t)},\,n\ge 0
\end{equation}
of quotients of generalized Fibonacci numbers $F_n(\t)$ defined by
\eqref{eq:genfib} are moment sequences of the measures $\mu_{\t}$ and
$\nu_{\t}=(2\sinh\t)\delta_1+\mu_{\t}$, where
\begin{equation}\label{eq:quofib1}
\mu_{\t}=e^{-\t}\delta_1 +2\cosh\t\sum_{k=1}^\infty
e^{-4k\t}\delta_{(-1)^ke^{-2k\t}}.
\end{equation}
In particular for $2\sinh\t=1$  the sequence $F_{n+1}/F_{n+2},n\ge
0$ is the moments of the probability measure
\begin{equation}\label{eq:quofib2}
\mu=\varphi\delta_1+\sqrt{5}\sum_{k=1}^\infty \varphi^{4k}\delta_{(-1)^k\varphi^{2k}},
\end{equation}
where $\varphi=(\sqrt{5}-1)/2$.
\end{cor}

\begin{rem} {\rm{Using the so-called Binet formula for the Fibonacci
      numbers, it is easy to see that $F_{n+1},n\ge 0$ is a moment
      sequence of the measure
$$
\tau=\frac{\sqrt{5}+1}{2\sqrt{5}}\delta_{(1+\sqrt{5})/2}+\frac{\sqrt{5}-1}{2\sqrt{5}}\delta_{(1-\sqrt{5})/2}.
$$
A similar formula holds for the  generalized Fibonacci numbers of
Ismail, see \cite[formula (2.2)]{Is1}.}}
\end{rem}

\begin{rem} {\rm{It was proved in \cite{Be2} that $F_\a/F_{n+\a},\,n\ge 0$ is the
moment sequence of a signed measure $\mu_\a$ with total mass 1. Here
$\a$ is a natural number and the signed measure
$\mu_\a$ is a probability measure precisely when $\a$ is an even number.
The orthogonal polynomials corresponding to $\mu_\a$
are little $q$-Jacobi polynomials, where
$q=(1-\sqrt{5})/(1+\sqrt{5})$. The results were used to prove
Richardson's formula for the elements in the inverse of the Filbert
matrix $(1/F_{1+i+j})$, cf. \cite{Ri}. These results were extended to
generalized Fibonacci numbers in \cite{Is1}. For an extension to quantum
integers see \cite{A:B}}}.
\end{rem}
\section{Proofs}

{\it Proof of Theorem~\ref{thm:1}.}
First of all, we find a closed expression for the denominators $(D_n)_n$ of the finite continued fractions $(s_n)_n$
defined in (\ref{partfrac}).

From (\ref{2suc}) and (\ref{2suc1}) we get
\begin{equation}\label{forn}
N_n=\frac{D_n-D_{n-2}}{a}, \quad n\ge 2.
\end{equation}
This recurrence can also be extended to $n=0, 1$ by defining $D_{-2}=1-aw$ and $D_{-1}=w$. Inserting (\ref{forn}) in
(\ref{2suc1}), we find that
\begin{equation}\label{dedn}
D_{n+2}=(2+ab)D_n-D_{n-2},\quad n\ge 0,
\end{equation}
with initial conditions $D_{-2}=1-aw$, $D_{-1}=w$, $D_{0}=1$ and $D_{1}=a+w$.

That means that the sequences $(D_{2n})_n$ and $(D_{2n+1})_n$ are both solutions of the difference equation
$$
x_{n+1}=(2+ab)x_n-x_{n-1},\quad n\ge 0,
$$
with initial conditions $x_{-1}=1-aw, w$, $x_{0}=1,a+w$, respectively. Any solution of this difference equation has the form $c_0q_0^n+c_1q_1^n$, where $q_0$ and $q_1$
are the solutions of $x^2-(2+ab)x+1=0$. We write
\begin{equation}\label{defq}
q=\frac{2+ab-\sqrt{a^2b^2+4ab}}{2},
\end{equation}
so that $q$ and $1/q$ are the solutions of $x^2-(2+ab)x+1=0$, and $0<q<1$. We then have that
there exist numbers $c_0,c_1,d_0,d_1$ such that
\begin{eqnarray}\label{exprdn}
D_{2n}&=&c_0q^{-n}+c_1q^n,\quad n\ge -1,\\
D_{2n+1}&=&d_0q^{-n}+d_1q^n,\quad n\ge -1.
\end{eqnarray}
Using the initial conditions $D_{-2}=1-aw$, $D_{-1}=w$, $D_{0}=1$ and $D_{1}=a+w$,
it is easy to calculate that
\begin{eqnarray}\label{exprcd}
c_0&=&\frac{1-q +qaw}{1-q^2}, \;\qquad c_1=\frac{q(1-q-aw)}{1-q^2} ,\\
d_0&=&\frac{a+(1-q)w}{1-q^2}, \qquad d_1=\frac{q((1-q)w-qa)}{1-q^2}.
\end{eqnarray}
Note that $c_0,d_0>0$.
Writing $\alpha=-c_1/c_0$ and $\beta =-d_1/(qd_0)$, we find
\begin{equation}\label{expralbe}
\alpha =\frac{q(aw-(1-q))}{qaw+1-q},\quad \beta =\frac{qa-(1-q)w}{a+(1-q)w},
\end{equation}
and it is clear that $|\alpha|,|\beta|<1$ because $a>0,0<q<1$ and
$w\ge 0$.

We need to establish some technical
properties of $\alpha $ and $\beta$, which we collect in

\begin{lem}\label{thm:lem2}
\begin{enumerate}
\item $\alpha \ge 0$ if and only if $w\ge -b/2+\sqrt{(b/2)^2+b/a}$.
\item Assume $w>0$. Then $-\beta \le \alpha$ if and only if $a\ge b$.
\item If $w=0$ then $\beta=-\alpha=q$.
\item $\alpha\ge \beta $ if and only if $w\ge -(a+b)/4+\sqrt{((a+b)/4)^2+1}$.
\end{enumerate}
\end{lem}

\begin{proof}

1. By \eqref{expralbe} we have that $\alpha\ge 0$ is equivalent to
$q+aw-1\ge 0$, hence to
$$ 1-aw\le q = \frac{2+ab-\sqrt{a^2b^2+4ab}}{2},
$$ 
or
$$
\sqrt{a^2b^2+4ab} \le 2aw+ab,
$$
which is equivalent to $aw^2+abw-b\ge 0$ because $a>0$. Since $a,b>0$ and $w\ge0$, we finally get that
$\alpha \ge 0$ if and only if 
$$
 w\ge -b/2+\sqrt{(b/2)^2+b/a}.
$$

2. According to \eqref{expralbe}, $-\beta \le \alpha$ if and only if
$$
\frac{(1-q)w-aq}{a+(1-q)w}\le \frac{q(aw-(1-q))}{qaw+1-q},
$$
and a straightforward computation gives that this is equivalent to
$$
w(1+q)(q^2-(2+a^2)q+1)\le 0.
$$
Using that $q$ satisfies $q^2-(2+ab)q+1=0$, the left-hand side of this inequality can be
reduced to $w(1+q)aq(b-a)$, and for $w>0$ the result follows.

3. Follows by inspection.
 
4. We similarly get that $\alpha\ge \beta$ if and only if
$$
w^2+w(a+b)/2-1\ge 0,
$$
which is equivalent to the given condition because $w\ge 0$.
\end{proof}

We now continue the proof of Theorem~\ref{thm:1}.

Using (\ref{forn}), we can write
$$
s_n=\frac{N_n}{D_n}= \frac{1}{a}\frac{D_n-D_{n-2}}{D_n}=\frac{1}{a}\left(1-\frac{D_{n-2}}{D_n}\right).
$$
From the formulas  (\ref{exprdn}) and (\ref{exprcd}), and taking into account
that $\alpha=-c_1/c_0$, we have
\begin{eqnarray*}
\frac{D_{2n-2}}{D_{2n}}&=&
\frac{c_0q^{-n+1}+c_1q^{n-1}}{c_0q^{-n}+c_1q^{n}}=q\frac{1-\alpha q^{2n-2}}{1-\alpha q^{2n}}\\
&=&q(1-\alpha q^{2n-2})\sum_{k=0}^\infty \alpha ^kq^{2nk}\\
&=&q\left(1+(1-q^{-2})\sum_{k=0}^\infty \alpha ^{k+1}q^{2n(k+1)}\right),
\end{eqnarray*}
showing that $s_{2n}=\int t^{2n}d\mu$, $n\ge 0$, where the measure $\mu $ is defined as
\begin{equation}\label{defmu}
\mu =\frac{1}{a}(1-q)\delta _1+\frac{1}{a}\left( \frac{1}{q}-q\right)\sum_{k=0}^\infty \alpha ^{k+1}\delta _{q^{k+1}}.
\end{equation}
In a similar way, it can be proved that $s_{2n+1}=\int t^{2n+1}d\nu$, $n\ge 0$, where the measure $\nu $ is defined as
\begin{equation}\label{defnu}
\nu =\frac{1}{a}(1-q)\delta _1+\frac{1}{a}\left( \frac{1}{q}-q\right)\sum_{k=0}^\infty \beta ^{k+1}\delta _{q^{k+1}}.
\end{equation}
Take now the reflected measures $\check \mu$ and $\check \nu$  of $\mu$ and $\nu$ with respect to the origin:
\begin{eqnarray*}
\check \mu &=&\frac{1}{a}(1-q)\delta _{-1}+\frac{1}{a}\left( \frac{1}{q}-q\right)\sum_{k=0}^\infty \alpha ^{k+1}\delta _{-q^{k+1}},\\
\check \nu &=&\frac{1}{a}(1-q)\delta _{-1}+\frac{1}{a}\left(
  \frac{1}{q}-q\right)\sum_{k=0}^\infty \beta ^{k+1}
\delta _{-q^{k+1}}.
\end{eqnarray*}
The measure $\rho =(\mu +\check \mu)/2 +(\nu -\check \nu)/2$ has then the same even moments as $\mu$ and the same
odd moments as $\nu$. That is, the $n$'th moment of $\rho$ is just $s_n$, $n\ge 0$. A simple computation shows that
the measure $\rho $ is given by (\ref{defro}).

It is clear that the measure $\rho $ is positive if and only if 
$$
\alpha ^{k+1}+\beta ^{k+1},\;\alpha ^{k+1}-\beta ^{k+1}\ge 0,\quad k\ge 0,
$$
i.e. if and only if
 $\alpha \ge 0$ and $-\alpha \le \beta\le \alpha $.
The second part of Theorem~\ref{thm:1} follows now by applying Lemma~\ref{thm:lem2}.

\hfill$\square$

{\it Proof of Corollary~\ref{thm:genfib}}. For $a=b=2\sinh\t>0$ and
$w=1/a$ we see that the three conditions of Theorem~\ref{thm:1} are
satisfied, so $s_n=s_n(a,a,1/a)$ is a moment sequence of a positive
measure $\rho=\mu_{\t}$. We prove that $s_n=F_{n+1}(\t)/F_{n+2}(\t)$
by induction. This formula holds for $n=0$ by inspection and clearly
$s_{n+1}=1/(a+s_n)$. We therefore find, assuming the formula for a fixed
$n$
$$
s_{n+1}=\frac{1}{a+ F_{n+1}(\t)/F_{n+2}(\t)}=\frac{F_{n+2}(\t)}{aF_{n+2}(\t)+F_{n+1}(\t)}=\frac{F_{n+2}(\t)}{F_{n+3}(\t)}.
$$
A small calculation shows that $\alpha=-\beta=q^2,q=e^{-2\t}$,
$e^{-\t}=(1-q)/a$ and $\mu_{\t}=\rho$ defined in \eqref{defro}
is given by \eqref{eq:quofib1}.
The formula for $\nu_{\t}$ follows easily from the recurrence equation
for $F_n(\t)$.
\hfill$\square$

\section{Concluding remarks}

Euler did not use 3-periodic continued fractions (nor any other periodicity bigger than 2)
to find rational approximations of square roots
of natural numbers: to use 3-periodic continued fractions is not more
useful than to use 1-periodic continued ones, and to use 4-periodic
continued fractions is not more
useful than to use 2-periodic ones,
 and so forth, cf.  \cite[p.322]{E}.

One can consider finite continued fractions like (\ref{partfrac}) for
a 3-periodic continued fraction
 of positive periods $a, b, c$.
Using the same approach as before, we can find three signed measures
$\mu_0$, $\mu_1$ and $\mu_2$ on $[-1,1]$ such that
the $(3n+i)$'th moment of $\mu_i$ is equal to $s_{3n+i}$, $i=0,1,2$,
$n\ge 0$. However, we are not able to construct  a measure on $\mathbb
R$
 from these three measures
having its $(3n+i)$'th moment equal to the $n$'th moment of $\mu_i$,
$i=0,1,2$. The same happens for $k$-periodicity when $k>2$.

By a result of Boas, cf. \cite[p. 138]{W}, see also \cite{Dur}, any 
sequence is a moment sequence
of a  signed measure. On the other hand we know that the Hankel
determinants of the moments of a positive measure are non-negative.
Computations indicate that if we consider the sequence
$(s_n(a,a,c,1))_n$, then some of the  Hankel determinants are negative 
except when $c=a$.

\end{document}